\documentclass[a4paper,reqno]{amsart}

\usepackage{amssymb}

\newtheorem{theorem}[equation]{Theorem}

\newtheorem{corollary}[equation]{Corollary}

\numberwithin{equation}{section}

\theoremstyle{definition}

\newtheorem*{example*}{Example}

\newtheorem*{remark*}{Remark}

\newcommand{\bZ}{{\mathbb Z}}

\newcommand{\frg}{{\mathfrak g}}

\newcommand{\frf}{{\mathfrak f}}
\newcommand{\fre}{{\mathfrak e}}

\newcommand{\calT}{{\mathcal T}}
\newcommand{\calTKK}{{\mathcal TKK}}

\newcommand{\subo}{_{\bar 0}}
\newcommand{\subuno}{_{\bar 1}}

\newcommand{\bil}{{\textup{b}}}

 \newcommand{\tri}{\mathfrak{tri}}

 \newcommand{\str}{\mathfrak{str}}
 \newcommand{\pstr}{\mathfrak{pstr}}
 
 \newcommand{\frsl}{{\mathfrak{sl}}}
 \newcommand{\frsp}{{\mathfrak{sp}}}
 \newcommand{\frso}{{\mathfrak{so}}}
 \newcommand{\frpsl}{{\mathfrak{psl}}}
 \newcommand{\frgl}{{\mathfrak{gl}}}
 \newcommand{\frpgl}{{\mathfrak{pgl}}}
 \newcommand{\frosp}{{\mathfrak{osp}}}

 \DeclareMathOperator{\der}{\mathfrak{der}}
 
 \DeclareMathOperator{\inder}{\mathfrak{inder}}
 
 \DeclareMathOperator{\End}{End}
 
 \DeclareMathOperator{\Mat}{Mat}

\def\bigstrut{\vrule height 12pt width 0ptdepth 2pt}

\def\hregleta{\hrule height .5pt}
\def\hreglon{\hrule height1pt}
\def\vreglon{\vrule height 12pt width1pt depth 4pt}

\def\hreglonfill{\leaders\hreglon\hfill}

\def\hregletafill{\leaders\hregleta\hfill}

\newenvironment{romanenumerate}
 {\begin{enumerate}
 
 }{\end{enumerate}}

\begin{document}

\title[The Supermagic Square and Jordan superalgebras]{The Supermagic Square in characteristic $3$ and Jordan superalgebras}

\author[Isabel Cunha]{Isabel Cunha$^{\diamond}$}
 \thanks{$^{\diamond}$ Supported by CMUC, Department of Mathematics, University of Coimbra}
 \address{Departamento de Matem\'atica, Universidade da Beira
 Interior,\newline 6200 Covilh\~a, Portugal}
 \email{icunha@mat.ubi.pt}

\author[Alberto Elduque]{Alberto Elduque$^{\star}$}
 \thanks{$^{\star}$ Supported by the Spanish Ministerio de
 Educaci\'{o}n y Ciencia
 and FEDER (MTM 2007-67884-C04-02) and by the
Diputaci\'on General de Arag\'on (Grupo de Investigaci\'on de
\'Algebra)}
 \address{Departamento de Matem\'aticas e
 Instituto Universitario de Matem\'aticas y Aplicaciones,
 Universidad de Zaragoza, 50009 Zaragoza, Spain}
 \email{elduque@unizar.es}

\dedicatory{Dedicated to Ivan Shestakov, on the occasion of his 60th birthday}

\date{February 22, 2008}



\begin{abstract}
Recently, the classical Freudenthal Magic Square has been extended over fields of characteristic $3$ with two more rows and columns filled with (mostly simple) Lie superalgebras specific of this characteristic. This \emph{Supermagic Square} will be reviewed and some of the simple Lie superalgebras that appear will be shown to be isomorphic to the Tits-Kantor-Koecher Lie superalgebras of some Jordan superalgebras.
\end{abstract}

\maketitle


\section*{Introduction}

The classical Freudenthal Magic Square, which contains in
characteristic $0$ the exceptional simple finite dimensional Lie
algebras, other than $G_2$, is usually constructed based on two
ingredients: a unital composition algebra and a central simple
degree $3$ Jordan algebra (see \cite[Chapter IV]{Schafer}). This
construction, due to Tits, does not work in characteristic $3$.

A more symmetric construction, based on two unital composition
algebras, which play symmetric roles, and their triality Lie
algebras, has been given recently by several authors (\cite{AF93},
\cite{BS1}, \cite{LM1} \cite{LM2}). Among other things, this construction has
the advantage of being valid too in characteristic $3$. Simpler
formulas for triality appear if symmetric composition algebras are
used, instead of the more classical unital composition algebras (\cite{EldIbero1,EldIbero2}).

But the characteristic $3$ presents an exceptional feature, as only
over fields of this characteristic there are nontrivial composition
superalgebras, which appear in dimensions $3$ and $6$. The unital such composition superalgebras were discovered by Shestakov \cite{She97}.
This fact
allows to extend Freudenthal Magic Square (\cite{CunEld1}) with the
addition of two further rows and columns, filled with (mostly
simple) Lie superalgebras, specific of characteristic $3$, which had appeared first (with one exception) in \cite{EldNew3} and \cite{EldModular}.

Most of the Lie superalgebras in characteristic $3$ that appear in the Supermagic Square have been shown to be related to degree three simple Jordan algebras in \cite{CunEld2}.

The aim of this paper is to show that some of the Lie superalgebras in the Supermagic Square are isomorphic to the Tits-Kantor-Koecher Lie superalgebras
of some distinguished  Jordan superalgebras.

More specifically, let $S_i$ denote the split para-Hurwitz algebra of dimension $i=1,2$ or $4$, and let $S$ be the para-Hurwitz superalgebra associated to the unital composition superalgebra $C$ (see Section 1 for definitions and notations). Let $\frg(S_i,S)$ be the corresponding entry in the Supermagic Square. Then $\frg(S_1,S)$ was shown in \cite{CunEld2} to be isomorphic to the Lie superalgebra of derivations of the Jordan superalgebra $J=H_3(C)$ of hermitian $3\times 3$ matrices over $C$. Here the following results will be proved:

\begin{romanenumerate}
\item
The Lie superalgebras $\frg(S_2,S)$ in the second row of the Supermagic Square will be shown to be isomorphic to the projective structure superalgebras of the Jordan superalgebras $J=H_3(C)$. Here the structure superalgebra is $\str(J)=L_J\oplus\der(J)$ and the projective structure superalgebra $\pstr(J)$ is the quotient of $\str(J)$ by its center. (See Theorem \ref{th:2ndrow} and Corollary \ref{co:2ndrowpstrJ}.)

\item
The Lie superalgebras $\frg(S_4,S)$ in the third row of the Supermagic Square will be shown to be isomorphic to the Tits-Kantor-Koecher Lie superalgebras of the Jordan superalgebras $J=H_3(C)$. (See Theorem \ref{th:Phi3} and Corollary \ref{co:3rdrowTKK}.)

\item
The Lie superalgebra $\frg(S_{1,2},S_{1,2})$ will be shown to be isomorphic to the Tits-Kantor-Koecher Lie superalgebra of the nine dimensional Kac Jordan superalgebra $K_9$. Note that the ten dimensional Kac Jordan superalgebra $K_{10}$ is no longer simple in characteristic $3$, but contains a nine dimensional simple ideal, which is $K_9$. (See Theorem \ref{th:gS12S12K9}.)

\item
The Lie superalgebra $\frg(S_1,S_{1,2})$ will be shown to be isomorphic to the Tits-Kantor-Koecher Lie superalgebra of the three dimensional Kaplansky superalgebra $K_3$. (See Corollary \ref{co:gS1S12K3}.)

\end{romanenumerate}

\smallskip

The paper is structured as follows. In Section 1 the construction of the Supermagic Square in terms of two symmetric composition superalgebras will be reviewed. Then the relationship of the Lie superalgebras in the first row of the Supermagic Square with the Lie superalgebras of derivations of the Jordan superalgebras $J=H_3(C)$ above, proven in \cite{CunEld2}, will be reviewed in Section 2. Section 3 will be devoted to the Lie superalgebras in the second and third rows of the Supermagic Square, while Section 4 will deal with the Lie superalgebra $\frg(S_{1,2},S_{1,2})$ and the nine dimensional Kac Jordan superalgebra $K_9$. It was Shestakov \cite{SheOberwolfach} who first noticed that $K_9$ is isomorphic to the tensor product (in the graded sense) of two copies of the three dimensional Kaplansky Jordan superalgebra $K_3$ (this was further developed in \cite{BenkartElduque}), and this is the key for the results in Section 4.

\smallskip

Unless otherwise stated, all the vector spaces and superspaces considered will be assumed to be finite dimensional over a ground field $k$ of characteristic $\ne 2$.

\bigskip
\section{The Supermagic Square}

A quadratic superform on a $\bZ_2$-graded vector space
$U=U\subo\oplus U\subuno$ over a field $k$ is a pair
$q=(q\subo,\bil)$ where $q\subo :U\subo\rightarrow k$ is a quadratic
form, and
 $\bil:U\times U\rightarrow k$ is a supersymmetric even bilinear form
such that $\bil\vert_{U\subo\times U\subo}$ is the polar form of $q\subo$:
\[
\bil(x\subo,y\subo)=q\subo(x\subo+y\subo)-q\subo(x\subo)-q\subo(y\subo)
\]
for any $x\subo,y\subo\in U\subo$.

The quadratic superform $q=(q\subo,\bil)$ is said to be
\emph{regular} if the bilinear form $\bil$  is
nondegenerate.

\smallskip

Then a superalgebra $C=C\subo\oplus C\subuno$ over $k$, endowed with
a regular quadratic superform $q=(q\subo,\bil)$, called the
\emph{norm}, is said to be a \emph{composition superalgebra} (see
\cite{EldOkuCompoSuper}) in case
\begin{subequations}\label{eq:norm}
\begin{align}
&q\subo(x\subo y\subo)=q\subo(x\subo)q\subo(y\subo),\label{eq:qcompo1}\\
&\bil(x\subo y,x\subo z)=q\subo(x\subo)\bil(y,z)=\bil(yx\subo,zx\subo),\label{eq:qcompo2}\\
&\bil(xy,zt)+(-1)^{  x   y  +
x
 z +  y   z }\bil(zy,xt)=(-1)^{
 y   z }\bil(x,z)\bil(y,t),\label{eq:qcompo3}
\end{align}
\end{subequations}
for any $x\subo,y\subo\in C\subo$ and homogeneous elements
$x,y,z,t\in C$. (As we are working in characteristic $\ne 2$, it is enough to consider equation \eqref{eq:qcompo3}.)

As usual, the expression $(-1)^{yz}$ equals $-1$ if the homogeneous elements $y$ and $z$ are both odd, otherwise, it equals $1$.

\smallskip

The unital composition superalgebras are termed \emph{Hurwitz
superalgebras}, while a composition superalgebra is said to be
\emph{symmetric} in case its bilinear form is associative, that is,
\[ \bil(xy,z)=\bil(x,yz),
\]
for any $x,y,z$.

\smallskip

Hurwitz algebras are the well-known algebras that generalize the
classical real division algebras of the real and complex numbers,
quaternions and octonions. Over any algebraically closed field $k$,
there are exactly four of them: $k$, $k\times k$, $\Mat_2(k)$ and
$C(k)$ (the split Cayley algebra), with dimensions $1$, $2$, $4$ and
$8$.

\smallskip

Only over fields of characteristic $3$  there appear nontrivial
Hurwitz superalgebras (see \cite{EldOkuCompoSuper}):

\begin{itemize}

\item Let $V$ be a two dimensional vector space over a field $k$,
endowed with a nonzero alternating bilinear form $\langle .\vert
.\rangle$ (that is $\langle v\vert v\rangle =0$ for any $v\in V$).  Consider the superspace $B(1,2)$ (see \cite{She97}) with
\begin{equation}\label{eq:B12a}
B(1,2)\subo =k1,\qquad\text{and}\qquad B(1,2)\subuno= V,
\end{equation}
endowed with the supercommutative multiplication given by
\[
1x=x1=x\qquad\text{and}\qquad uv=\langle u\vert v\rangle 1
\]
for any $x\in B(1,2)$ and $u,v\in V$, and with the quadratic
superform $q=(q\subo,\bil)$ given by:
\begin{equation}\label{eq:B12b}
q\subo(1)=1,\quad \bil(u,v)=\langle u\vert v\rangle,
\end{equation}
for any $u,v\in V$. If the characteristic of $k$ is equal to $3$, then
$B(1,2)$ is a Hurwitz superalgebra (\cite[Proposition
2.7]{EldOkuCompoSuper}).

\smallskip

\item Moreover, with $V$ as before, let $f\mapsto \bar f$ be the
associated symplectic involution on $\End_k(V)$ (so $\langle
f(u)\vert v\rangle =\langle u\vert\bar f(v)\rangle$ for any $u,v\in
V$ and $f\in\End_k(V)$). Consider the superspace $B(4,2)$ (see
\cite{She97}) with
\begin{equation}\label{eq:B42}
B(4,2)\subo=\End_k(V),\qquad\text{and}\qquad B(4,2)\subuno=V,
\end{equation}
with multiplication given by the usual one (composition of maps) in
$\End_k(V)$, and by
\[
\begin{split}
&v\cdot f=f(v)=\bar f\cdot v \in V,\\
&u\cdot v=\langle .\vert u\rangle v\in \End_k(V)
\end{split}
\]
for any $f\in\End_k(V)$ and $u,v\in V$, where $\langle .\vert u\rangle v$ denotes the endomorphism $w\mapsto \langle w\vert u\rangle v$; and with quadratic superform
$q=(q\subo,\bil)$ such that
\[
q\subo(f)=\det(f),\qquad\bil(u,v)=\langle u\vert v\rangle,
\]
for any $f\in \End_k(V)$ and $u,v\in V$. If the
characteristic is equal to $3$, $B(4,2)$ is a Hurwitz superalgebra
(\cite[Proposition 2.7]{EldOkuCompoSuper}).

\end{itemize}

\smallskip

Given any Hurwitz superalgebra $C$ with norm $q=(q\subo,\bil)$, its
standard involution is given by
\[
x\mapsto \bar x=\bil(x,1)1-x.
\]
A new product can be defined on $C$ by means of
\begin{equation}\label{eq:paraHurwitz}
x\bullet y=\bar x\bar y.
\end{equation}
The resulting superalgebra, denoted by $\bar C$, is called the
\emph{para-Hurwitz superalgebra} attached to $C$, and it turns out to be a symmetric composition superalgebra.

\smallskip

Given a symmetric composition superalgebra $S$, its \emph{triality
Lie superalgebra} $\tri(S)=\tri(S)\subo\oplus\tri(S)\subuno$ is
defined by:
\begin{multline*}
\tri(S)_{\bar i}=\{ (d_0,d_1,d_2)\in\frosp(S,q)^3_{\bar i}:\\
d_0(x\bullet y)=d_1(x)\bullet y+(-1)^{i  x }x\bullet
d_2(y)\ \forall x,y\in S\subo\cup S\subuno\},
\end{multline*}
where $\bar i= \bar 0,\bar 1$, and $\frosp(S,q)$ denotes the
associated orthosymplectic Lie superalgebra. The bracket in
$\tri(S)$ is given componentwise.

Now, given two symmetric composition superalgebras $S$ and $S'$, one can form (see
\cite[\S 3]{CunEld1}, or \cite{EldIbero1} for the non-super situation) the Lie superalgebra:
\begin{equation}\label{eq:gSS'}
\frg=\frg(S,S')=\bigl(\tri(S)\oplus\tri(S')\bigr)\oplus\bigl(\oplus_{i=0}^2
\iota_i(S\otimes S')\bigr),
\end{equation}
where $\iota_i(S\otimes S')$ is just a copy of $S\otimes S'$
($i=0,1,2$),  with bracket given by:

\smallskip
\begin{itemize}
\item the Lie bracket in $\tri(S)\oplus\tri(S')$, which thus becomes  a Lie subsuperalgebra of $\frg$,
\smallskip

\item $[(d_0,d_1,d_2),\iota_i(x\otimes
 x')]=\iota_i\bigl(d_i(x)\otimes x'\bigr)$,
\smallskip

\item
 $[(d_0',d_1',d_2'),\iota_i(x\otimes
 x')]=(-1)^{  d_i'   x }\iota_i\bigl(x\otimes d_i'(x')\bigr)$,
\smallskip

\item $[\iota_i(x\otimes x'),\iota_{i+1}(y\otimes y')]=(-1)^{
x'   y }
 \iota_{i+2}\bigl((x\bullet y)\otimes (x'\bullet y')\bigr)$ (indices modulo
 $3$),
\smallskip

\item $[\iota_i(x\otimes x'),\iota_i(y\otimes y')]=
 (-1)^{  x   x' +  x   y'  +
   y   y' }
 \bil'(x',y')\theta^i(t_{x,y})$ \newline \null\hspace{2.5 in} $+
 (-1)^{  y   x' }
 \bil(x,y)\theta'^i(t'_{x',y'})$,

\end{itemize}
\smallskip

\noindent
for any $i=0,1,2$ and homogeneous $x,y\in S$, $x',y'\in S'$,
$(d_0,d_1,d_2)\in\tri(S)$, and $(d_0',d_1',d_2')\in\tri(S')$. Here
$\theta$ denotes the natural automorphism
$\theta:(d_0,d_1,d_2)\mapsto (d_2,d_0,d_1)$ in $\tri(S)$, while $t_{x,y}$ is defined by
\begin{equation}\label{eq:txy}
t_{x,y}=\bigl(\sigma_{x,y},\tfrac{1}{2}\bil(x,y)1-r_xl_y,\tfrac{1}{2}\bil(x,y)1-l_xr_y\bigr)
\end{equation}
with $l_x(y)=x\bullet y$, $r_x(y)=(-1)^{xy}y\bullet x$, and
\begin{equation}\label{eq:sigmaxy}
\sigma_{x,y}(z)=(-1)^{yz}\bil(x,z)y-(-1)^{x(y+z)}\bil(y,z)x
\end{equation}
for homogeneous $x,y,z\in S$. Also $\theta'$
and $t'_{x',y'}$ denote the analogous elements for $\tri(S')$.

\smallskip

Over a field $k$ of characteristic $3$, let $S_r$ ($r=1$, $2$, $4$ or $8$) denote the para-Hurwitz algebra attached to the split Hurwitz
algebra of dimension $r$ (this latter algebra being either $k$, $k\times k$,
$\Mat_2(k)$ or $C(k)$). Also, denote by
$S_{1,2}$ the para-Hurwitz superalgebra
$\overline{B(1,2)}$, and by $S_{4,2}$ the
para-Hurwitz superalgebra $\overline{B(4,2)}$. Then the Lie superalgebras $\frg(S,S')$, where $S,S'$ run over $\{S_1,S_2,S_4,S_8,S_{1,2},S_{4,2}\}$, appear in Table \ref{ta:supermagicsquare}, which has been obtained in \cite{CunEld1}.

\begin{table}[h!]
$$
\vbox{\offinterlineskip
 \halign{\hfil\ $#$\ \hfil&%
 \vreglon #%
 &\hfil\ $#$\ \hfil&\hfil\ $#$\ \hfil
 &\hfil\ $#$\ \hfil&\hfil\ $#$\ \hfil&%
 \vrule  depth 4pt width .5pt #%
 &\hfil\ $#$\ \hfil&\hfil\ $#$\ \hfil\cr
 \bigstrut &width 0pt&S_1&S_2&S_4&S_8&\omit%
    \vrule height 8pt depth 4pt width .5pt&S_{1,2}&S_{4,2}\cr
 &\multispan8{\hreglonfill}\cr
 S_1&&\frsl_2&\frpgl_3&\frsp_6&\frf_4&&\frpsl_{2,2}&\frsp_6\oplus (14)\cr
 \bigstrut S_2&& &\omit$\frpgl_3\oplus \frpgl_3$&\frpgl_6&\tilde \fre_6&%
     &\bigl(\frpgl_3\oplus\frsl_2\bigr)\oplus\bigl(\frpsl_3\otimes (2)\bigr)&
    \frpgl_6\oplus (20)\cr
 \bigstrut S_4&& & &\frso_{12}&\fre_7&
   &\bigl(\frsp_6\oplus\frsl_2\bigr)\oplus\bigl((13)\otimes (2)\bigr)
    &\frso_{12}\oplus spin_{12}\cr
 \bigstrut S_8&& & & &\fre_8&
    &\bigl(\frf_4\oplus\frsl_2\bigr)\oplus\bigl((25)\otimes (2)\bigr)&
      \fre_7\oplus (56)\cr
 \multispan9{\hregletafill}\cr
 \bigstrut S_{1,2}&& & & & & & \frso_7\oplus 2spin_7 &\frsp_8\oplus(40)\cr
 \bigstrut S_{4,2}&& & & & & & & \frso_{13}\oplus spin_{13}\cr}}
$$
\bigskip
\caption{Supermagic Square (characteristic
$3$)}\label{ta:supermagicsquare}
\end{table}

Since the construction of $\frg(S,S')$ is symmetric, only the
entries above the diagonal are needed. In Table
\ref{ta:supermagicsquare}, $\frf_4,\fre_6,\fre_7,\fre_8$ denote the
simple exceptional classical Lie algebras, $\tilde\fre_6$ denotes a
$78$ dimensional Lie algebras whose derived Lie algebra is the $77$
dimensional simple Lie algebra $\fre_6$ in characteristic
$3$. The even and odd parts of the nontrivial superalgebras in the
table which have no counterpart in the classification in
characteristic $0$ (\cite{Kac-Lie}) are displayed, $spin$ denotes the
spin module for the corresponding orthogonal Lie algebra, while
$(n)$ denotes a module of dimension $n$, whose precise description is given in \cite{CunEld1}. Thus, for example,
$\frg(S_4,S_{1,2})$ is a Lie superalgebra whose even part is
(isomorphic to) the direct sum of the symplectic Lie algebra
$\frsp_6$ and of $\frsl_2$, while its odd part is the tensor
product of a $13$ dimensional module for $\frsp_6$ and the
natural $2$ dimensional module for $\frsl_2$.

A precise description of these modules and of the Lie superalgebras
as Lie superalgebras with a Cartan matrix is given in \cite{CunEld1}. All the
inequivalent Cartan matrices for these simple Lie superalgebras are listed in \cite{BGL}.

With the exception of $\frg(S_{1,2},S_{4,2})$, all these superalgebras have appeared previously in \cite{EldNew3} and \cite{EldModular}.

\bigskip
\section{Jordan superalgebras}\label{se:Jordan}

Given any Hurwitz superalgebra $C$ over our ground field $k$, with norm $q=(q\subo,\bil)$ and standard
involution $x\mapsto \bar x$, the superalgebra $H_3(C)$ of
$3\times 3$ hermitian matrices over $C$, under the superinvolution given by $(a_{ij})^*=(\bar a_{ji})$, is a Jordan superalgebra under the symmetrized product
\begin{equation}\label{eq:Jproduct}
x\circ y= \frac{1}{2}\bigl( xy+(-1)^{xy}yx\bigr).
\end{equation}

Let us consider the associated para-Hurwitz superalgebra $S=\bar C$,
with multiplication $a\bullet b=\bar a\bar b$ for any $a,b\in C$.
Then,
\begin{equation}\label{eq:JH3C}
\begin{split}
J=H_3(C)&=\left\{ \begin{pmatrix} \alpha_0 &\bar a_2& a_1\\
  a_2&\alpha_1&\bar a_0\\ \bar a_1&a_0&\alpha_2\end{pmatrix} :
  \alpha_0,\alpha_1,\alpha_2\in k,\ a_0,a_1,a_2\in S\right\}\\[6pt]
 &=\bigl(\oplus_{i=0}^2 ke_i\bigr)\oplus
     \bigl(\oplus_{i=0}^2\iota_i(S)\bigr),
\end{split}
\end{equation}
where
\begin{equation}\label{eq:eisiotas}
\begin{aligned}
e_0&= \begin{pmatrix} 1&0&0\\ 0&0&0\\ 0&0&0\end{pmatrix}, &
 e_1&=\begin{pmatrix} 0&0&0\\ 0&1&0\\ 0&0&0\end{pmatrix}, &
 e_2&= \begin{pmatrix} 0&0&0\\ 0&0&0\\ 0&0&1\end{pmatrix}, \\
 \iota_0(a)&=2\begin{pmatrix} 0&0&0\\ 0&0&\bar a\\
 0&a&0\end{pmatrix},&
 \iota_1(a)&=2\begin{pmatrix} 0&0&a\\ 0&0&0\\
 \bar a&0&0\end{pmatrix},&
 \iota_2(a)&=2\begin{pmatrix} 0&\bar a&0\\ a&0&0\\
 0&0&0\end{pmatrix},
\end{aligned}
\end{equation}
for any $a\in S$. Identify $ke_0\oplus ke_1\oplus ke_2$ to $k^3$ by
means of $\alpha_0e_0+\alpha_1e_1+\alpha_2e_2\simeq
(\alpha_0,\alpha_1,\alpha_2)$. Then the supercommutative
multiplication \eqref{eq:Jproduct} becomes:
\begin{equation}\label{eq:Jniceproduct}
\left\{\begin{aligned}
 &(\alpha_0,\alpha_1,\alpha_2)\circ(\beta_1,\beta_2,\beta_3)=
    (\alpha_0\beta_0,\alpha_1\beta_1,\alpha_2\beta_2),\\
 &(\alpha_0,\alpha_1,\alpha_2)\circ \iota_i(a)
  =\frac{1}{2}(\alpha_{i+1}+\alpha_{i+2})\iota_i(a),\\
 &\iota_i(a)\circ\iota_{i+1}(b)=\iota_{i+2}(a\bullet b),\\
 &\iota_i(a)\circ\iota_i(b)=2\bil(a,b)\bigl(e_{i+1}+e_{i+2}\bigr),
\end{aligned}\right.
\end{equation}
for any $\alpha_i,\beta_i\in k$, $a,b\in S$, $i=0,1,2$, and where
indices are taken modulo $3$.

\smallskip

In \cite{CunEld2} it is  shown that the Lie superalgebra of
derivations of $J$ is naturally isomorphic to the Lie superalgebra
$\frg(S_1,S)$ in the first row of the Supermagic Square.

This is well-known for algebras, as $\frg(S_1,S)$ is isomorphic to
the Lie algebra $\calT(k,H_3(C))$ obtained by means of Tits
construction (see \cite{EldIbero1} and \cite{BS2}), and this latter
algebra is, by its own construction, the derivation algebra of
$H_3(C)$. What was done in \cite[Section 3]{CunEld2} is to make explicit
this isomorphism $\frg(S_1,S)\cong \der J$ and extend it to
superalgebras.

To begin with, \eqref{eq:Jniceproduct} shows that $J$ is graded over
$\bZ_2\times \bZ_2$ with:
\[
J_{(0,0)}=k^3,\quad J_{(1,0)}=\iota_0(S),\quad
 J_{(0,1)}=\iota_1(S),\quad J_{(1,1)}=\iota_2(S)
\]
and, therefore, $\der J$ is accordingly graded over
$\bZ_2\times\bZ_2$:
\[
(\der J)_{(i,j)}=\{ d\in\der J: d\bigl(J_{(r,s)}\bigr)\subseteq
J_{(i+r,j+s)}\ \forall r,s=0,1\}.
\]
Moreover, the zero component is (\cite[Lemmas 3.4 and 3.5]{CunEld2}):
\[
(\der J)_{(0,0)}=\{ d\in \der J : d(e_i)=0\ \forall
i=0,1,2\},
\]
and the linear map given by
\[
\begin{split}
\tri(S)&\longrightarrow (\der J)_{(0,0)}\\
 (d_0,d_1,d_2)&\mapsto D_{(d_0,d_1,d_2)},
\end{split}
\]
where
\begin{equation}\label{eq:Dd0d1d2}
\left\{\begin{aligned} &D_{(d_0,d_1,d_2)}(e_i)=0,\\
   &D_{(d_0,d_1,d_2)}\bigl(\iota_i(a)\bigr)=
   \iota_i\bigl(d_i(a)\bigr)
   \end{aligned}\right.
\end{equation}
for any $i=0,1,2$ and $a\in S$, is an isomorphism.

\smallskip

Given any two elements $x,y$ in a Jordan superalgebra, the commutator (in the graded sense) of the left multiplications by $x$ and $y$:
\begin{equation}\label{eq:dxy}
d_{x,y}=[L_x,L_y]
\end{equation}
is a derivation. These derivations are called inner derivations.
For any $i=0,1,2$ and $a\in S$, consider the following inner
derivation of the Jordan superalgebra $J$:
\begin{equation}\label{eq:Dia}
D_i(a)=2d_{\iota_i(a),e_{i+1}}=2\bigl[ L_{\iota_i(a)},L_{e_{i+1}}\bigr]
\end{equation}
(indices modulo $3$), where $L_x$ denotes the multiplication by $x$
in $J$. Note that the restriction of $L_{e_i}$ to
$\iota_{i+1}(S)\oplus\iota_{i+2}(S)$ is half the identity, so the
inner derivation $\bigl[ L_{\iota_i(a)},L_{e_{i}}\bigr]$ is trivial
on $\iota_{i+1}(S)\oplus\iota_{i+2}(S)$, which generates $J$. Hence
\begin{equation}\label{eq:Liotaiaei}
\bigl[ L_{\iota_i(a)},L_{e_{i}}\bigr]=0
\end{equation}
for any $i=0,1,2$ and $a\in S$. Also, $L_{e_0+e_1+e_2}$ is the
identity map, so the bracket $\bigl[ L_{\iota_i(a)},L_{e_0+e_{1}+e_2}\bigr]$ is $0$
and hence
\begin{equation}\label{eq:Liotaiaei2}
D_i(a)=2\bigl[ L_{\iota_i(a)},L_{e_{i+1}}\bigr]=-2\bigl[
L_{\iota_i(a)},L_{e_{i+2}}\bigr].
\end{equation}

A straightforward computation with \eqref{eq:Jniceproduct} gives
\begin{equation}\label{eq:Diaaction}
\begin{split}
&D_i(a)(e_i)=0,\ D_i(a)(e_{i+1})=\frac{1}{2} \iota_i(a),\
  D_i(a)(e_{i+2})=-\frac{1}{2}\iota_i(a),\\
&D_i(a)\bigl(\iota_{i+1}(b)\bigr)=-\iota_{i+2}(a\bullet b),\\
&D_i(a)\bigl(\iota_{i+2}(b)\bigr)=(-1)^{\lvert a\rvert\lvert
  b\rvert}\iota_{i+1}(b\bullet a),\\
&D_i(a)\bigl(\iota_i(b)\bigr)=2\bil(a,b)(-e_{i+1}+e_{i+2}),
\end{split}
\end{equation}
for any $i=0,1,2$ and any homogeneous elements $a,b\in S$.

Denote by $D_i(S)$ the linear span of the $D_i(a)$'s, $a\in S$. Then the remaining components of the $\bZ_2\times\bZ_2$-grading of $\der J$ are given by (\cite[Lemma 3.11]{CunEld2}):
\[
(\der J)_{(1,0)}=D_0(S),\quad (\der J)_{(0,1)}=D_1(S),\quad
(\der J)_{(1,1)}=D_2(S).
\]

Therefore, the $\bZ_2\times\bZ_2$-grading of $\der J$ becomes
\begin{equation}\label{eq:derJDS}
\der J=D_{\tri(S)}\oplus\bigl(\oplus_{i=0}^2 D_i(S)\bigr)
\end{equation}

On the other hand, $S_1=k1$, with $1\bullet 1=1$ and $\bil(1,1)=2$,
so $\tri(S_1)=0$ and for the para-Hurwitz superalgebra $S$:
\[
\begin{split}
\frg(S_1,S)&=\tri(S)\oplus\bigl(\oplus_{i=0}^2 \iota_i(S_1\otimes
      S)\bigr)\\
   &=\tri(S)\oplus\bigl(\oplus_{i=0}^2\iota_i(1\otimes S)\bigr).
\end{split}
\]

\begin{theorem}\label{th:gS1SderJ} (See \cite[Theorem 3.13]{CunEld2}) Let $S$ be a para-Hurwitz
superalgebra over $k$ and let $J$ be
the Jordan superalgebra of $3\times 3$ hermitian matrices over the
associated Hurwitz superalgebra. Then the linear map:
\begin{equation}\label{eq:PhigS1SderJ}
\Phi:\frg(S_1,S)\longrightarrow \der J,
\end{equation}
such that
\[
\begin{split}
\Phi\bigl((d_0,d_1,d_2)\bigr)&=D_{(d_0,d_1,d_2)},\\
\Phi\bigl(\iota_i(1\otimes a)\bigr)&=D_i(a),
\end{split}
\]
for any $i=0,1,2$, $a\in S$ and $(d_0,d_1,d_2)\in\tri(S)$, is an
isomorphism of Lie superalgebras.
\end{theorem}

\smallskip

The Lie superalgebra $d_{J,J}=[L_J,L_J]$ is the Lie superalgebra $\inder J$
of inner derivations of $J$. It turns out that $(\der
J)_{(r,s)}=(\inder J)_{(r,s)}$ for $(r,s)\ne (0,0)$, while
\[
(\inder
J)_{(0,0)}=\sum_{i=0}^2\bigl[L_{\iota_i(S)},L_{\iota_i(S)}\bigr]=
D_{\sum_{i=0}^2\theta^i(t_{S,S})}=\Phi\bigl(\sum_{i=0}^2\theta^i(t_{S,S})\bigr)
\]
(recall that $\theta\bigl((d_0,d_1,d_2)\bigr)=(d_2,d_0,d_1)$ for any
$(d_0,d_1,d_2)\in\tri(S)$).

In characteristic $3$, $\tri(S)=\sum_{i=0}^2\theta^i(t_{S,S})$ if
$\dim S=1,4$ or $8$ (\cite{EldIbero1}), and the same happens with
$\tri(S_{1,2})$ and $\tri(S_{4,2})$, because of \cite[Corollaries
2.12 and 2.23]{CunEld1}, while for $\dim S=2$, $\tri(S)$ has
dimension $2$ and $\sum_{i=0}^2\theta^i(t_{S,S})=t_{S,S}$ has
dimension $1$ (see \cite{EldIbero1}). In characteristic $\ne 3$,
$\tri(S)=\sum_{i=0}^2\theta^i(t_{S,S})$ always holds.

\begin{corollary}
(See \cite[Corollary 3.15]{CunEld2})
Let $S$ be a para-Hurwitz (super)algebra over $k$, and let $J$  be the Jordan (super)algebra of
$3\times 3$ hermitian matrices over the associated Hurwitz
(super)algebra. Then $\der J$ is a simple Lie (super)algebra that
coincides with $\inder J$ unless the characteristic is $3$ and $\dim
S=2$. In this latter case $\inder J$ coincides with $[\der J,\der
J]$, which is a codimension $1$ simple ideal of $\der J$.
\end{corollary}

\bigskip
\section{The second and third rows in the Supermagic Square}

The second row of the Supermagic Square is formed by the Lie superalgebras $\frg(S,S')$, where $S$ and $S'$ are symmetric composition superalgebras with $\dim S=2$. Let $S_2$ be the split two dimensional para-Hurwitz algebra. Thus, $S_2$ is the para-Hurwitz algebra attached to the unital composition algebra $K=k\times k$, whose standard involution is given by $\overline{(\alpha,\beta)}=(\beta,\alpha)$ for any $\alpha,\beta\in k$. Then with $1=(1,1)$ and $u=(1,-1)$, the multiplication and norm in $S_2$ are given by:
\[
\begin{aligned}
1\bullet 1&=1,& 1\bullet u&=u\bullet 1=-u,& u\bullet u&=1,\\
q(1)&=1,& \bil(1,u)&=\bil(u,1)=0,& q(u)&=-1.
\end{aligned}
\]
Moreover, the triality Lie algebra of $S_2$ is (see \cite[Corollary 3.4]{EldIbero1}) the Lie algebra:
\[
\tri(S_2)=\{(\alpha_0\sigma,\alpha_1\sigma,\alpha_2\sigma): \alpha_0,\alpha_1,\alpha_2\in k\ \text{and}\ \alpha_0+\alpha_1+\alpha_2=0\},
\]
where the linear map
\begin{equation}\label{eq:sigma}
\sigma=\frac{1}{2}\sigma_{1,u}
\end{equation}
(recall \eqref{eq:sigmaxy}) satisfies $\sigma(1)=u$ and $\sigma(u)=1$. Note that the element $t_{1,u}\in\tri(S_2)$ defined in \eqref{eq:txy} is given by:
\begin{equation}\label{eq:t1u}
t_{1,u}=(2\sigma,-\sigma,-\sigma).
\end{equation}

Let now $J$ be the Jordan superalgebra considered in equation \eqref{eq:JH3C}. Its \emph{structure Lie superalgebra} (or Lie multiplication superalgebra) $\str J$ is the subalgebra of the general Lie superalgebra $\frgl(J)$ spanned by the Lie superalgebra of derivations $\der J$ of $J$ and by the space $L_J$ of left multiplications by elements in $J$. Since $[L_J,L_J]$ is a subalgebra of $\der J$, it follows that
\[
\str J=\der J\oplus L_J.
\]
Then the center of this Lie superalgebra is spanned by $L_1$ (the identity map). We will consider too the projective structure Lie superalgebra $\pstr J$, which is defined as the quotient of $\str J$ modulo its center:
\[
\pstr J=\str J/kL_1 =\str J/kI,
\]
where $I$ denotes the identity map on $J$.

Then with the notations introduced in Section \ref{se:Jordan}, we have:

\begin{theorem}\label{th:2ndrow}
Let $S$ be a para-Hurwitz superalgebra over  $k$ and let $J$ be the Jordan superalgebra of $3\times 3$ hermitian matrices over the associated Hurwitz superalgebra. Then the isomorphism $\Phi$ in equation \eqref{eq:PhigS1SderJ} extends to the following isomorphism of Lie superalgebras
\begin{equation}\label{eq:Phi2}
\Phi_2:\frg(S_2,S)\longrightarrow \pstr J
\end{equation}
where:
\begin{itemize}
\item for any $(d_0,d_1,d_2)\in\tri(S)$, $\Phi_2\bigl((d_0,d_1,d_2)\bigr)= D_{(d_0,d_1,d_2)} + kI$,
\item for any $\alpha_0,\alpha_1,\alpha_2\in k$ with $\alpha_0+\alpha_1+\alpha_2=0$, the image under $\Phi$ of the element
    $(\alpha_0\sigma,\alpha_1\sigma,\alpha_2\sigma)\in\tri(S_2)$ is
   $L_{\alpha_2e_1-\alpha_1e_2}+kI\,
   \bigl(=L_{\alpha_1e_0-\alpha_0e_1}+kI=L_{\alpha_0e_2-\alpha_2e_1}+kI\bigr)$,
\item and for any $i=0,1,2$, $a\in S$, $\Phi_2\bigl(\iota_i\bigl((\alpha 1+\beta u)\otimes a)\bigr)=\bigl(\alpha D_i(a)+\beta L_{\iota_i(a)}\bigr)+kI$.
\end{itemize}
\end{theorem}
\begin{proof}
The proof is obtained by straightforward computations using the results in Section \ref{se:Jordan}. Thus, for instance, for any $\alpha_i\in k$ ($i=0,1,2$) with $\alpha_0+\alpha_1+\alpha_2=0$, and any $\alpha,\beta\in k$ and $a\in S$, we get:
\[
\bigl[(\alpha_0\sigma,\alpha_1\sigma,\alpha_2\sigma),\iota_0\bigl((\alpha 1+\beta u)\otimes a\bigr)\bigr]=\alpha_0\iota_0\bigl((\beta 1+\alpha u)\otimes a\bigr),
\]
which maps under $\Phi_2$ to
\[
\alpha_0\bigl(\beta D_0(a)+\alpha L_{\iota_0(a)}\bigr)+kI.
\]
On the other hand, we have:
\[
\begin{split}
\bigl[\Phi_2\bigl((\alpha_0\sigma,\alpha_1\sigma,\alpha_2\sigma))&,\Phi_2\bigl(\iota_0\bigl((\alpha 1+\beta u)\otimes a\bigr)\bigr)\bigr]\\
 &=\bigl[L_{\alpha_2 e_1-\alpha_1 e_2},\alpha D_0(a)+\beta L_{\iota_0(a)}\bigr]+kI\\
 &=\left(-\alpha[D_0(a),L_{\alpha_1e_0-\alpha_0e_1}]
   -\beta[L_{\iota_0(a)},L_{\alpha_1e_0-\alpha_0e_1}]\right)+kI\\
 &=\alpha_0\left(\alpha L_{D_0(a)(e_1)}+\beta D_0(a)\right)+kI\\
 &=\alpha_0\left(\beta D_0(a)+\alpha L_{\iota_0(a)}\right)+kI,
\end{split}
\]
where equations \eqref{eq:Dia}, \eqref{eq:Liotaiaei}, \eqref{eq:Liotaiaei2} and \eqref{eq:Diaaction} have been used.

In a similar vein, for any homogeneous $a,b\in S$, the element
\[
[\iota_0(1\otimes a),\iota_0(u\otimes b)]=\bil(a,b)t_{1,u}=\bil(a,b)(2\sigma,-\sigma,\sigma)
\]
(recall equations \eqref{eq:sigma} and \eqref{eq:t1u}) maps under $\Phi_2$ to $\bil(a,b)L_{e_2-e_1}+kI$, while we have
\[
\begin{split}
\bigl[\Phi_2\bigl(\iota_0(1\otimes a)\bigr)&,
        \Phi_2\bigl(\iota_0(u\otimes b)\bigr)\bigr]\\
 &=[D_{\iota_0(a)},L_{\iota_0(b)}]+kI\\
 &=L_{D_0(a)(\iota_0(b))}+kI=\bil(a,b)L_{-e_1+e_2}+kI,
\end{split}
\]
as required.

The remaining computations needed to prove that $\Phi_2$ is an isomorphism are similar to the ones above, and will be omitted.
\end{proof}

\begin{corollary}\label{co:2ndrowpstrJ}
The Lie superalgebras $\frg(S_2,S_{1,2})$ and $\frg(S_2,S_{4,2})$ in the Supermagic Square in characteristic $3$ are isomorphic, respectively, to the projective structure Lie superalgebras of the Jordan superalgebras of hermitian $3\times 3$ matrices over the unital composition superalgebras $B(1,2)$ and $B(4,2)$.
\end{corollary}

\medskip

Let us turn now our attention to the third row of the Supermagic Square.

Thus, let $Q$ be a quaternion algebra (that is, a four dimensional unital composition algebra) over $k$, with multiplication denoted by juxtaposition, and denote by $\bar Q$ the para-Hurwitz algebra with multiplication given by $x\bullet y=\overline{xy}=\bar y \bar x$. Note that $\bar Q$ is the para-Hurwitz algebra attached to the opposite algebra of $Q$ (which is isomorphic to $Q$). In case $Q$ is split, then it is isomorphic to the algebra of $2\times 2$ matrices $\Mat_2(k)$.

According to \cite[Corollary 3.4]{EldIbero1}, the triality Lie algebra of $\bar Q$ splits as:
\begin{equation}\label{eq:pi012}
\tri(\bar Q)=\ker \pi_0\oplus\ker\pi_1\oplus\ker\pi_2,
\end{equation}
where $\pi_i:\tri(\bar Q)\rightarrow \frso(Q)$ is the projection onto the $i$th component. Moreover, let $Q^0$ denote the subspace of zero trace elements in $Q$, that is, the subspace orthogonal to the unity element. Then (\cite[Corollary 3.4]{EldIbero1}) we have:
\[
\ker \pi_0=\{(0,l_a\tau,-r_a\tau): a\in Q^0\},
\]
where $l_a$ and $r_a$ denote the left and right multiplications in $\bar Q$ and $\tau:x\mapsto \bar x$ is the standard involution of $Q$. Therefore, for any $a\in Q^0$ and $x\in Q$, since $\bar a=\tau(a)=-a$, we get:
\[
\begin{split}
l_a\tau(x)&=a\bullet\bar x=\overline{a\bar x}=x\bar a=-xa=-R_a(x),\\
r_a\tau(x)&=\bar x\bullet a=\overline{\bar x a}=\bar ax=-ax=-L_a(x),
\end{split}
\]
where $L_a$ and $R_a$ denote the left and right multiplications by $a$ in $Q$. Hence the ideal $\ker\pi_0$ above becomes:
\[
\ker \pi_0=\{ (0,-R_a,L_a):a\in Q^0\}.
\]
and, similarly, $\ker\pi_1=\{(L_a,0,-R_a):a\in Q^0\}$ and $\ker\pi_2=\{(-R_a,L_a,0):a\in Q^0\}$.

Now, for any $a,b,x\in Q$, $q(x)=x\bar x=\bar x x$, so $\bil(a,b)=a\bar b+b\bar a=\bar a b+b\bar a$ and hence we have:
\[
\begin{split}
&r_al_b(x)=\overline{\overline{bx}a}=\overline{\bar x\bar b a}=\bar a b x=L_{\bar a b}(x),\\
&l_ar_b(x)=\overline{a\overline{xb}}=\overline{a\bar b\bar x}=x b\bar a =R_{b\bar a}(x),\\
&\sigma_{a,b}(x)=\bil(a,x)b-\bil(b,x)a=(a\bar x+x\bar a)b-a(\bar b x+\bar x b)
    = (-L_{a\bar b}+R_{\bar a b})(x),\\
&\sigma_{a,b}(x)=\bil(a,x)b-\bil(b,x)a=b(\bar a x+\bar x a)-(b\bar x+x\bar b)a
    = (L_{b\bar a}-R_{\bar b a})(x),\\
&\frac{1}{2}\bil(a,b)x-r_al_b(x)
  =\frac{1}{2}(\bar a b+\bar ba)x-L_{\bar a b}(x)
   =\frac{1}{2}L_{\bar b a-\bar a b}(x),\\
&\frac{1}{2}\bil(a,b)x-l_ar_b(x)
   =\frac{1}{2} x(a\bar b+b\bar a)-R_{b\bar a}(x)
   =\frac{1}{2}R_{a\bar b-b\bar a}(x).
\end{split}
\]
Therefore,
the element $t_{a,b}\in \tri(\bar Q)$ in \eqref{eq:txy} becomes:
\[
\begin{split}
t_{a,b}&=\bigl(\sigma_{a,b},\frac{1}{2}\bil(a,b)1-r_al_b,\frac{1}{2}\bil(a,b)-l_ar_b\bigr)\\
 &=\frac{1}{2}\bigl( -L_{a\bar b-b\bar a}+R_{\bar ab-\bar ba},
    -L_{\bar ab-\bar ba},R_{a\bar b-b\bar a}\bigr).
\end{split}
\]

It must be noticed that the subspace $Q^0$ of trace zero elements is a three dimensional simple Lie algebra under the commutator, and that any three dimensional simple Lie algebra appears in this way. Moreover, given any Jordan algebra $H$, Tits considered in \cite{Tits62} the Lie algebra defined on the vector space
\begin{equation}\label{eq:TQH}
\calT(Q,H)=\bigl(Q^0\otimes H\bigr)\oplus \der H
\end{equation}
endowed with the bracket given by:
\[
\begin{split}
\bullet&\ \textrm{the restriction to $\der H$ is the commutator in $\der H$,}\\
\bullet&\ [d,a\otimes x]=a\otimes d(x),\\
\bullet&\ [a\otimes x,b\otimes y]=\bigl([a,b]\otimes
xy\bigr)-2\bil(a,b)d_{x,y},
\end{split}
\]
for any $a,b\in Q^0$, $x,y\in H$, and $d\in\der H$. (Recall the definition of $d_{x,y}$ in \eqref{eq:dxy}.)

In the split case
$Q=\Mat_2(k)$,  the resulting Lie algebra is the well-known
Tits-Kantor-Koecher Lie algebra $\calTKK(H)$ of the Jordan algebra $H$. Besides, all the arguments involved work in the super setting, and thus $H$ can be taken to be any Jordan superalgebra.

\smallskip

Let us return our attention to the Jordan superalgebra $J$ of hermitian $3\times 3$ matrices over a unital composition superalgebra $C$, with associated para-Hurwitz superalgebra denoted by $S$, as in equation \eqref{eq:JH3C}. Consider the Lie superalgebra $\frg(\bar Q, S)$ in equation \eqref{eq:gSS'}. As a vector space, $\bar Q$ splits into the direct sum $\bar Q=k1\oplus Q^0$, and this gives the following decomposition of $\frg(\bar Q,S)$:
\[
\frg(\bar Q,S)=\frg(S_1,S)\oplus \tri(\bar Q)\oplus\bigl(\oplus_{i=0}^2\iota_i(Q^0\otimes S)\bigr)
\]
Also recall the decomposition in equation \eqref{eq:pi012}.

Then, as in Theorem \ref{th:2ndrow}, the isomorphism $\Phi$ in equation \eqref{eq:PhigS1SderJ} can be extended to $\frg(\bar Q,S)$:

\begin{theorem}\label{th:Phi3}
Let $S$ be a para-Hurwitz superalgebra over  $k$ and let $J$ be the Jordan superalgebra of $3\times 3$ hermitian matrices over the associated Hurwitz superalgebra. Then the isomorphism $\Phi$ in equation \eqref{eq:PhigS1SderJ} extends to the following isomorphism of Lie superalgebras
\begin{equation}\label{eq:Phi3}
\Phi_3:\frg(\bar Q,S)\longrightarrow \calT(Q,J)=\bigl(Q^0\otimes J)\oplus\der J,
\end{equation}
where:
\begin{itemize}
\item the restriction of $\Phi_3$ to $\frg(S_1,S)$ coincides with $\Phi$ in Theorem \ref{th:gS1SderJ},
\item for any $a\in Q^0$, the elements $(0,-R_a,L_a)\in \ker\pi_0$, $(L_a,0,-R_a)\in\ker\pi_1$ and $(-R_a,L_a,0)\in \ker\pi_2$ map, respectively, to $\frac{1}{2}(a\otimes e_0)$, $\frac{1}{2}(a\otimes e_1)$, and $\frac{1}{2}(a\otimes e_2)$, (recall that equation \eqref{eq:pi012} shows that $\tri(\bar Q)=\ker\pi_0\oplus\ker\pi_1\oplus\ker\pi_2$,)
\item for any $i=0,1,2$, $a\in Q^0$ and $x\in S$, $\Phi_3\bigl(\iota_i(a\otimes x)= -\frac{1}{2}a\otimes \iota_i(x)$.
\end{itemize}
\end{theorem}

The proof is obtained by straightforward computations and thus will be omitted.

\bigskip

Since $\calT(Q,J)$ is the Tits-Kantor-Koecher Lie superalgebra of the Jordan superalgebra $J$ in case $Q$ is the split quaternion algebra, the next corollary follows at once:

\begin{corollary}\label{co:3rdrowTKK}
The Lie superalgebras $\frg(S_4,S_{1,2})$ and $\frg(S_4,S_{4,2})$ in the Supermagic Square in characteristic $3$ are isomorphic, respectively, to the Tits-Kantor-Koecher Lie superalgebras of the Jordan superalgebras of hermitian $3\times 3$ matrices over the unital composition superalgebras $B(1,2)$ and $B(4,2)$.
\end{corollary}

\bigskip
\section{The Lie superalgebra $\frg(S_{1,2},S_{1,2})$}

The tiny \emph{Kaplansky superalgebra} (\cite{Kap,McC}) is the three dimensional Jordan superalgebra $K_3=K\subo\oplus K\subuno$, with $K\subo=ke$ and $K\subuno=kx+ky$, and with multiplication given by:
\begin{gather*}
e^2=e,\quad ex=\frac{1}{2}x=xe,\quad ey=\frac{1}{2}y=ye,\\
xy=e=-yx,\quad x^2=0=y^2.
\end{gather*}

On the other hand, the simple ten dimensional Kac superalgebra $K_{10}$ was originally constructed in \cite{Kac-Jordan} over algebraically closed fields of characteristic $0$ by Lie-theoretical methods from a $3$-grading of the exceptional Lie superalgebra $F(4)$. In characteristic $3$, $K_{10}$ is no longer simple but possesses a simple ideal $K_9$ of dimension $9$. Shestakov (\cite{SheOberwolfach} unpublished) noticed that $K_9$ is isomorphic to the tensor product (as superalgebras) $K_3\otimes K_3$. Later on, it was proven in \cite{BenkartElduque} that $K_{10}$ appears as a direct sum $k1\oplus (K_3\otimes K_3)$, with a natural multiplication, in any characteristic, and in particular, if the characteristic is $3$, then $K_9=K_3\otimes K_3$.

\smallskip

Assume in the remaining of this section that the characteristic of our ground field $k$ is $3$.

\smallskip

Take the unital composition superalgebra $B(1,2)$ in \eqref{eq:B12a} and its para-Hurwitz counterpart $S_{1,2}$. Then take a symplectic basis $\{u,v\}$ of $(S_{1,2})\subuno=V$, so that $\{1,u,v\}$ is a basis of $S_{1,2}$. Since the characteristic is $3$, the multiplication of $S_{1,2}$ is given by (see \eqref{eq:paraHurwitz}):
\begin{gather*}
1\bullet 1=\bar 1\bar 1=1,\quad 1\bullet x=\bar 1\bar x=1(-x)=-x=\frac{1}{2}x=x\bullet 1\quad \textrm{for any $x\in (S_{1,2})\subuno$},\\
u\bullet v=\bar u\bar v=(-u)(-v)=1,\quad u\bullet u=0=v\bullet v,
\end{gather*}
and therefore $S_{1,2}$ is just the tiny Kaplansky superalgebra $K_3$.

Hence, Kac superalgebra $K_9$ can be identified with the superalgebra $S_{1,2}\otimes S_{1,2}$.

Now, the triality Lie superalgebra of $S_{1,2}$ is computed in \cite[Theorem 5.6]{EldOkuCompoSuper} (see also \cite[Corollary 2.12]{CunEld1}):
\begin{equation}\label{eq:triS12}
\tri(S_{1,2})=\{(d,d,d): d\in \frosp(S_{1,2},\bil)\}.
\end{equation}
That is, $\tri(S_{1,2)}$ is isomorphic to the orthosymplectic Lie superalgebra on the three dimensional vector superspace $S_{1,2}$ relative to the polar form of its norm. Also, in \cite[Theorem 5.8]{EldOkuCompoSuper} it is proven that the Lie superalgebra of derivations of $S_{1,2}$ is the whole orthosymplectic Lie superalgebra $\frosp(S_{1,2},\bil)$.

In \cite[Theorem 2.8]{BenkartElduque} it is proven that the Lie superalgebra of derivations of the Kac superalgebra $K_9=K_3\otimes K_3$ is the direct sum of the Lie superalgebras of derivations of the two copies of the tiny Kaplansky superalgebra involved. That is, we have:
\[
\der K_9=(\der K_3\otimes I)\oplus(I\otimes \der K_3),
\]
where $I$ denotes the identity map and given any homogeneous linear maps $\varphi,\psi\in\frgl(K_3)$, $\varphi\otimes \psi$ is the linear endomorphism
in $\frgl(K_9)=\frgl(K_3\otimes K_3)$ given by:
\[
(\varphi\otimes\psi)(x\otimes y)=(-1)^{\psi x}(\varphi(x)\otimes\psi(y))
\]
for any homogeneous elements $x,y\in K_3$.

\smallskip

Consider now the quaternion algebra $Q$ over $k$ with a basis $\{1,e_0,e_1,e_2\}$, where $1$ is the unity element, and with
\[
e_i^2=-1,\qquad e_ie_{i+1}=-e_{i+2}=-e_{i+1}e_i,
\]
(indices modulo $3$). Its norm is the regular quadratic form $q$ with $q(1)=1=q(e_i)$, for $i=0,1,2$ and where the basis above is orthogonal.
Since the characteristic is $3$, $q(e_0+e_1+e_2)=0$, so the norm of $Q$ represents $0$ and hence $Q$ is the split quaternion algebra, that is, it is isomorphic to $\Mat_2(k)$. (An explicit isomorphism can be easily constructed.)

To avoid confusion, let us denote by $\bil_q$ the polar form of $q$.
The subspace of zero trace elements is $Q^0=ke_0+ke_1+ke_2$.

The Lie superalgebra $\calT(Q,K_9)$ (see equation \eqref{eq:TQH}), which is isomorphic to the Tits-Kantor-Koecher Lie superalgebra of $K_9$ as $Q$ is split, is given by
\[
\calT(Q,J)=\bigl(Q^0\otimes K_9\bigr)\oplus \der K_9,
\]
and hence it decomposes as:
\[
\calT(Q,K_9)=\bigl(\oplus_{i=0}^2 e_i\otimes (K_3\otimes K_3)\bigr)\oplus\bigl((\der K_3\otimes I)\oplus(I\otimes \der K_3)\bigr),
\]
that is, the direct sum of three copies of the tensor product $K_3\otimes K_3$ (or $S_{1,2}\otimes S_{1,2}$) and two copies of $\der K_3$, which is isomorphic to the orthoysimplectic Lie superalgebra $\frosp(S_{1,2},\bil)$, exactly the situation that occurs for the Lie superalgebra $\frg(S_{1,2},S_{1,2})=\tri(S_{1,2})\oplus \tri(S_{1,2})\oplus\bigl(\oplus_{i=0}^2\iota_i(S_{1,2}\otimes S_{1,2})$.

\begin{theorem}\label{th:gS12S12K9}
Let $k$ be a field of characteristic $3$. Then the Lie superalgebra $\frg(S_{1,2},S_{1,2})$ is isomorphic to the Tits-Kantor-Koecher Lie superalgebra of the Kac superalgebra $K_9$.
\end{theorem}
\begin{proof}
It has been checked above that both Lie superalgebras, $\frg(S_{1,2},S_{1,2})$ and the Tits-Kantor-Koecher Lie superalgebra $\calT(Q,K_9)$, split as vector superspaces into direct sums of isomorphic summands. Let us consider the explicit linear isomorphism:
\[
\Psi:\frg(S_{1,2},S_{1,2})\longrightarrow \calT(Q,K_9),
\]
given by:
\begin{itemize}
\item $\Psi\bigl((d,d,d)\bigr)=d\otimes I\in\der K_9$, for any $(d,d,d)$ in the first copy of $\tri(S_{1,2})$ in $\frg(S_{1,2},S_{1,2})$,
\item $\Psi(\bigl((d',d',d')\bigr)=I\otimes d'$, for any $(d',d',d')$ in the second copy of $\tri(S_{1,2})$ in $\frg(S_{1,2},S_{1,2})$,
\item $\Psi\bigl(\iota_i(x\otimes x')=e_i\otimes (x\otimes x')$, for any $i=0,1,2$ and $x,x'\in S_{1,2}$,
\end{itemize}
and let us check that it is an isomorphism of Lie superalgebras.

In order to prove that $\Psi$ is indeed a Lie superalgebra isomorphism, the only nontrivial point is to prove that
\[
\Psi\bigl([\iota_i(x\otimes x'),\iota_i(y\otimes y')]\bigr)=
 [\Psi(\iota_i(x\otimes x')),\Psi(\iota_i(y\otimes y'))]
\]
for any homogeneous, $x,x',y,y'\in S_{1,2}=K_3$ and $i=0,1,2$. The symmetry of the constructions shows that it is enough to deal with $i=0$. But the description of $\tri(S_{1,2})$ in \eqref{eq:triS12}, together with equations \eqref{eq:gSS'}, \eqref{eq:txy} and \eqref{eq:sigmaxy} give:
\[
\begin{split}
[\iota_0(x\otimes x'),\iota_0(y\otimes y')]
 &=(-1)^{xx'+xy'+yy'}\bil(x',y')t_{x,y}+(-1)^{x'y}\bil(x,y)t_{x',y'}\\
 &=(-1)^{xx'+xy'+yy'}\bil(x',y')(\sigma_{x,y},\sigma_{x,y},\sigma_{x,y})\\
 &\qquad\qquad
 +(-1)^{x'y}\bil(x,y)(\sigma_{x',y'},\sigma_{x',y'},\sigma_{x',y'}),
\end{split}
\]
which maps, under $\Psi$ to:
\[
\begin{split}
\Psi\Bigl([\iota_0(x\otimes x')&,\iota_0(y\otimes y')]\Bigr)\\
 &=(-1)^{xx'+xy'+yy'}\bil(x',y')(\sigma_{x,y}\otimes I)
 +(-1)^{x'y}\bil(x,y)(I\otimes \sigma_{x',y'}),\\
 &=(-1)^{x'y}\left((\sigma_{x,y}\otimes \bil(x',y')I)+
    (\bil(x,y)I\otimes \sigma_{x',y'})\right)\in\der K_9.
\end{split}
\]
On the other hand, we have:
\[
\begin{split}
[\Psi\bigl(\iota_0(x\otimes x')\bigr)&,\Psi\bigl(\iota_0(y\otimes y')\bigr)]\\
 &=[e_0\otimes (x\otimes x'),e_0\otimes(y\otimes y')]\\
 &=-2b_q(e_0,e_0)[L_{x\otimes x'},L_{y\otimes y'}]\quad\textrm{(recall the product in \eqref{eq:TQH})}\\
 &=-[L_{x\otimes x'},L_{y\otimes y'}]\quad\textrm{(as $q(e_0)=1$, so $\bil_q(e_0,e_0)=2=-1$)}\\
 &=-(-1)^{x'y}\frac{1}{2}\left(([L_x,L_y]\otimes \bil(x',y')I)+(\bil(x,y)I\otimes [L_{x'},L_{y'}])\right)\\
 &\hspace{120pt}\textrm{(by \cite[(2.3)]{BenkartElduque})}\\
 &=(-1)^{x'y}\left((\sigma_{x,y}\otimes \bil(x',y')I)+(\bil(x,y)I\otimes \sigma_{x',y'})\right)\\
 &\hspace{120pt} \textrm{(because of \cite[(1.6)]{BenkartElduque})}\\
 &=\Psi\Bigl([\iota_0(x\otimes x'),\iota_0(y\otimes y')]\Bigr),
\end{split}
\]
as required.
\end{proof}

\smallskip

The Lie superalgebra $\frg(S_1,S_{1,2})$ is a subalgebra of $\frg(S_{1,2},S_{1,2})$. The restriction of the isomorphism $\Psi$ in the proof of Theorem \ref{th:gS12S12K9} gives our last result:

\begin{corollary}\label{co:gS1S12K3}
Let $k$ be a field of characteristic $3$. Then the Lie superalgebra $\frg(S_1,S_{1,2})$ in the Supermagic Square is isomorphic to the Tits-Kantor-Koecher Lie superalgebra of the tiny Kaplansky superalgebra $K_3$.
\end{corollary}

Note that the Lie superalgebra $\frg(S_1,S_{1,2})$ is known to be isomorphic to the projective special Lie superalgebra $\frpsl_{2,2}$ (see \cite{CunEld1}).


\providecommand{\bysame}{\leavevmode\hbox
to3em{\hrulefill}\thinspace}

\end{document}